\documentclass[12pt,singlespace,oneside]{article}
\usepackage{amsfonts}

\usepackage{amsmath}
\usepackage{latexsym}

\begin{document}

\title{Real Version of Calculus of Complex Variable (I): Weierstrass Point of
View\thanks{%
Dedicated to Professor Zbigniew Oziewicz, whom I owe much for educating me
on some sophisticated mathematical issues.} \ }
\author{Jose G. Vargas\thanks{%
138 Promontory Rd., Columbia, SC 29209, USA. josegvargas@earthlink.net} \ }
\date{}
\maketitle

\begin{abstract}
A very small amount of K\"{a}hler algebra (i.e. Clifford algebra of
differential forms) in the real plane makes $x+ydxdy$ emerge as a factor
between the differentials of the Cartesian and polar coordinates, largely
replacing the concept of complex variable. The integration on closed curves
of closed 1-forms on multiply connected regions takes us directly to a real
plane version of the theorem of residues. One need not resort to anything
like differentiation and integration with respect to $x+ydxdy.$ It is a
matter of algebra and integration of periodic functions. We then derive
Cauchy's integral formulas, including the ones for the derivatives.
Additional complex variable theory of general interest for phyicists are
then trivial.

The approach is consistent with the Wierstrass point of view: power series
expansions, even if explicit expressions are not needed. By design, this
approach cannot replace integrations that yield complex results. These can
be obtained with an approach based on the Cauchy point of view, where the
Cauchy-Riemann conditions come first and the theorem of residues comes last
(Paper to follow).
\end{abstract}

\section{Introduction}

In this paper, we reach the theorem of residues on multiply connected
domains of the real plane, thus bringing Cauchy's theory to the fold of real
analysis. It results from a simple and direct continuation of a corollary to
Stokes theorem in not simply connected regions generated by the removal of
poles of integrands. We assume that these are isolated ones.

The key element in this as in a follow-up paper is the replacement of $z$
with the differential form $x+ydxdy.$ Readers who know if credit has to be
given to others for developing similar work (\textit{without
differentiations with respect to a Clifford variable!}), please contact this
author.

In this Weierstrass approach, we get the theorem of residues as if the
Cauchy calculus had never existed. We then go into developing Cauchy's
special integral formulas. We use the term special since it does not apply
to (the real form) of imaginary integrals \ The more general framework in
next paper's approach, which we shall call Cauchy's, will allow us to obtain
those ``complex results''. The perspective of those approaches is due to
Henri Cartan \cite{HCartan}, but for the fact that our algebra is real and
no new concept of differentiation is needed.

This paper can be used for a two hour course. In the first hour, one would
derive the theorem of residues and pose a set of exercise for practice. In a
second hour, one would resolve problems that students might have had solving
those exercises and would present Cauchy's integral formulas, with
applications. But this should not distract us from the conceptual progress
made.

The words in the language of integration are differential forms. For the use
that most physicists make of the theory of complex variable, it is a
historical accident that happened at a time when differential forms were not
known (much less their Clifford algebras). The genius of Cauchy allowed him
to get powerful results without the right tools to address them.

\section{Extension of the real calculus in the plane}

\subsection{Essence of the theorem of residues}

In the real plane, Stokes theorem allows one to replace arbitrary
integration contours (defined here as closed curves, not as boundaries) with
circles centered at isolated singularities. In addition, we can make the
radius go towards zero without changing the value of the integral. Hence the
problem of integration becomes one of computing limits of integrals that
depend only on the angular coordinate.

Let $\alpha $ be a closed scalar-valued differential 1-form, i.e. $d\alpha
=0 $. This is equivalent to the condition of zero curl of a vector field
with the same components as $\alpha $. In polar coordinates, $\alpha $ is
written as%
\begin{equation}
\alpha =h(\rho ,\phi )d\rho +j(\rho ,\phi )d\phi .
\end{equation}%
We wish to integrate $\alpha $ on curves that enclose only isolated poles
but go through none. We may replace the integral with a sum of non
overlapping circles with the same orientation, each centered at and
containing only one pole. When limits exist, the fact that $\rho $ remains
constant on centered circles (indicated with the subscript $0$, and omitting
a subscript to label the poles) allows us to write that integral with a sum
of integrals over the small circles,%
\begin{equation}
\oint \alpha =\sum \lim_{\rho _{0}\rightarrow 0}\oint j(\rho _{0}\text{,}%
\phi )d\phi .
\end{equation}%
The idea behind the theorem of residues is that integration of the Fourier
power series of $j(\rho _{0}$,$\phi )$ over $2\pi $ implies%
\begin{equation}
\oint \alpha =2\pi \sum \lim_{\rho _{0}\rightarrow 0}a_{0},
\end{equation}%
where $a_{0}$ represents the constant term in the Fourier series of the
expansion of $j(\rho _{0}$,$\phi )$ around each pole. But its computation
involves performing the integral that one intends to solve. Clifford algebra
allows us to easily convert this trigonometric problem into one involving
Cartesian coordinates. Its solution then lies in obtaining a limit. That is
in a nutshell the theorem of residues.

\subsection{Theorem of residues}

The K\"{a}hler algebra (i.e. Clifford algebra of differential forms) is
defined by%
\begin{equation}
dx^{i}dx^{j}+dx^{j}dx^{i}=2g^{ij},
\end{equation}%
where $g^{ij}$ is the metric \cite{K62}. In Cartesian and polar coordinates
in any real plane:%
\begin{eqnarray}
(dx)^{2} &=&(dy)^{2}=(d\rho )^{2}=1,\text{ \ \ \ \ \ }(d\phi )^{2}=\frac{1}{%
\rho ^{2}} \\
dxdy &=&-dydx,\text{ \ }\ \ \ \ \ \ \ \ \ \ (dxdy)^{2}=-1.
\end{eqnarray}%
The complex looking inhomogeneous differential form%
\begin{equation}
z\equiv x+ydxdy,
\end{equation}%
emerges from the relation between ($d\rho ,d\phi $) and ($dx,dy$):%
\begin{eqnarray}
d\phi &=&\frac{xdy-ydx}{x^{2}+y^{2}}=\frac{x-ydxdy}{x^{2}+y^{2}}dy=\frac{1}{%
x+ydxdy}dy=\frac{1}{z}dy, \\
d\rho &=&\frac{xdx+ydy}{(x^{2}+y^{2})^{1/2}}=\rho \frac{x-ydxdy}{x^{2}+y^{2}}%
dx=\frac{\rho }{x+ydxdy}dx=\frac{\rho }{z}dx.
\end{eqnarray}%
By virtue of the second equation (6), it is clear that%
\begin{equation}
z^{\pm m}=(x+ydxdy)^{\pm m}=\rho ^{\pm m}e^{m\phi dxdy}=\rho ^{\pm m}(\cos
m\phi \pm dxdy\sin m\phi ),
\end{equation}%
for integer $m.$

Like $z$ itself, functions $F(z)$ take values in the even subalgebra, whose
elements are of the form%
\begin{equation}
u(x,y)+v(x,y)dxdy,
\end{equation}%
This subalgebra is commutative. In the full algebra, we have%
\begin{equation}
(u+vdxdy)\alpha =\alpha (u+vdxdy)^{\ast },\text{ \ \ \ \ \ }(u+vdxdy)^{\ast
}\equiv u-vdxdy,
\end{equation}%
and, in particular,%
\begin{equation}
z^{\ast }=x-ydxdy,\text{ \ \ \ }z^{\ast }=\rho ^{2}/z
\end{equation}

Given two differential 1-forms $\alpha $ and $\beta $, define%
\begin{equation}
\alpha \cdot \beta \equiv \frac{1}{2}(\alpha \beta +\beta \alpha ).
\end{equation}%
One readily obtains%
\begin{equation}
d\rho \cdot d\phi =0,\text{ \ \ }d\phi \cdot d\phi =\frac{1}{\rho ^{2}},
\end{equation}%
which yields%
\begin{equation}
j=\rho ^{2}(\alpha \cdot d\phi ).
\end{equation}

Solving for $x$ and $y$ in the system of equations (7) and (13) and
replacing $\rho $ with $\rho _{0}$ for integration on circles centered at
the origin, we have 
\begin{equation}
x=\frac{z+z^{\ast }}{2}=\frac{1}{2}\left( z+\frac{\rho _{0}^{2}}{z}\right)
,\ \ \ \ \ \ \ \ \ \ \ y=\frac{1}{2dxdy}\left( z-\frac{\rho _{0}^{2}}{z}%
\right) .
\end{equation}%
On those circles, functions of $x$ and $y$ become functions of $z$\ ---not
also $z^{\ast }$. Trigonometric functions reduce to the case just considered
since $\cos \phi =x/\rho _{0}$ and $\sin \phi =y/\rho _{0}$. We can always
write the pull-back $f(x,y)dx+g(x,y)dy$ of $\alpha $ to Cartesian
coordinates as%
\begin{equation}
k(x,y)dx+g(x,y)dy=wdx,\text{ \ \ \ \ \ }w\equiv k-gdxdy.
\end{equation}%
We proceed to compute $j$: 
\begin{eqnarray}
j &=&\rho ^{2}(wdx)\cdot (\frac{1}{z}dy)=\frac{\rho ^{2}}{2}\left[ wdx\frac{1%
}{z}dy+\frac{1}{z}dywdx\right] =  \notag \\
&=&\frac{\rho ^{2}}{2}\left[ w\frac{1}{z^{\ast }}dxdy+\frac{1}{z}w^{\ast
}dydx\right] =\frac{1}{2}\left[ wzdxdy+w^{\ast }z^{\ast }(dxdy)^{\ast }%
\right] =  \notag \\
&=&(wzdxdy)^{(0)}=-(wz)^{(2)},
\end{eqnarray}%
where the superscripts $``0"$ and $``2"$ stand here for $u$ and $v$ in $wz$
(not in $w$ here). They play the role of real and imaginary parts of $wz$.
For circles centered at pole, $z_{0}$, the role of $z$ will be played by%
\begin{equation}
z^{\prime }\equiv z-z_{0}.
\end{equation}

We shall use the term analytic to refer to functions $f$ of $z$ given by a
power series, an example being the sine function, where \textit{sine }is an
abbreviation for the series. We shall be interested in meromorphic functions 
$f$. They are defined as being analytic in the open set obtained by removing
isolated points from the set of definition of $f$. Such is the case with
integer power expansions extended by a finite number of negative power
terms, like in the quotient of the sine function by a polynomial. The zeroes
of the polynomial are called the poles of the function. In any quotient of
functions, the zeroes of the denominator are poles (i.e. points of
divergence) of the quotient function, except possibly if the numerator also
has a zero at the same point.

Assume that $w$ admits or directly is a series of integer powers (positive,
zero and negative) of $z^{\prime }$ (Fractional powers are not periodic over
a circle!). If circling just a pole of first order, we have%
\begin{equation}
\oint \alpha =\lim_{\rho _{0}\rightarrow 0}\oint j(\rho _{0}\text{,}\phi
)d\phi =-2\pi \lim_{z\rightarrow z_{0}}[(z-z_{0})w]_{c}^{(2)}],
\end{equation}%
the subscript $c$ referring to the constant term in the series for $%
[(z-z_{0})w]^{(2)}$.

For poles of arbitrary (integer) order, $m$, that limit does not pick the
constant coefficient of the series for $(z-z_{0})w$.\ We shall later
consider the standard way of picking that term in the Cauchy calculus. It is
worth recalling that, for most integrals of interest, we can proceed with
the known alternative that follows.

For a given function $\mathfrak{F}(z),$ let $m$ be the order (i.e. the
smallest positive integer such that $\lim_{z\rightarrow z_{i}}\mathfrak{F}%
\cdot (z-z_{i})^{m}$ is a number $a_{m}$ different from zero). If 
\begin{equation}
\mathfrak{F}-\frac{a_{m}}{z^{m}}
\end{equation}%
still has a pole, it is of order $p$ not higher than $m-1.$ We then compute $%
a_{p}$, and so on, until we reach the constant term. It may be zero in
particular.

\subsection{On applying the theorem of residues}

In dealing with integrals on the whole real line, the usual first step in
the application of the theorem of residues is the obvious equality%
\begin{equation}
\int_{-\infty }^{+\infty }H(x)dx=\int_{-\infty }^{+\infty }wdx,
\end{equation}%
where $w$ is defined as $H(z).$ If $H(z)$ goes sufficiently fast to zero on
a semicircle $\Gamma $ at the infinity of the upper or lower half of the
real plane, we use it to make with the $x$ axis a closed curve.

Consider, for instance, the integral%
\begin{equation}
\int_{-\infty }^{+\infty }\text{ }\frac{1}{\left( x^{2}+1\right) ^{2}}dx.
\end{equation}%
The poles are $z=\pm dxdy$, of order $m=2$, and $1/(z^{2}+1)^{2}$ is clearly
meromorphic. We shall later compute it through Cauchy's integral formula,
but, if we want to use the theorem of residues, we compute $a_{2}$ at $%
z=dxdy $:%
\begin{equation}
a_{2}=\lim_{z\rightarrow dxdy}\frac{(z-dxdy)^{2}}{\left( z^{2}+1\right) ^{2}}%
=\lim_{z\rightarrow dxdy}\frac{(z-dxdy)^{2}}{(z-dxdy)^{2}(z+dxdy)^{2}}=-%
\frac{1}{4}.
\end{equation}%
$\left( \mathfrak{F}-\frac{a_{2}}{z^{\prime 2}}\right) $ will now have a
pole of first order at $z=dxdy$, i.e. at the point ($0,1$). In order to
minimize clutter when computing $z^{\prime }\left( \mathfrak{F}-\frac{a_{2}}{%
z^{\prime 2}}\right) $, we exceptionally use the symbol $i$ as abbreviation
for $dxdy.$ Hence%
\begin{eqnarray}
z^{\prime }\left( \mathfrak{F}-\frac{a_{2}}{z^{\prime 2}}\right) &=&(z-i)%
\left[ \frac{1}{(z^{2}+1)^{2}}+\frac{1}{4(z-i)^{2}}\right] =  \notag \\
&=&\frac{z^{3}+iz^{2}+5z-3i}{4(z^{2}+1)^{2}}=\frac{z+3i}{4(z+i)^{2}}
\end{eqnarray}%
The factor of $i$ in the limit\ $z\rightarrow i$ of this expression is $%
-1/4. $ Multiplying by $-2\pi $, we obtain the value $\pi /2$ for the
integral (24).

Up to this point in this section, we have assumed the integer powers
expansion. It is a trivial matter to show that these expansions satisfy the
Cauchy Riemann conditions,%
\begin{equation}
u,_{x}=v,_{y}\text{ \ \ \ \ \ \ \ \ \ \ \ \ \ \ \ \ \ \ }u,_{y}=-v,_{x}
\end{equation}%
which are independent of each other.

Suppose that, rather than generating $w$ ($=H(z)$) by replacing $x$ with $z$
in $H(x)$, we originally had a closed $1-$form $\alpha $ that we were to
integrate on a closed curve in the plane. For $\alpha $ to be closed, we
only need $k,_{y}=g,_{x}$. We rewrite (18) as 
\begin{equation}
wdx\equiv k(x,y)dx+g(x,y)dy=(k-gdxdy)dx.
\end{equation}%
The identifications%
\begin{equation}
u\equiv k,\text{ \ \ \ \ \ \ \ \ \ }v\equiv -g
\end{equation}%
follow. Condition $k,_{y}=g,_{x}$ is only the second of conditions (27). The
method of residues is not justified for all closed differential forms but
only for those that also satisfy $k,_{x}=-g,_{y}$, since our argument
applied only to functions that are integer power series, even if not
expressed as such (remember that, for example, \textit{sine} is an
abbreviation for a particular series).

\section{Integral formulas}

\subsection{Cauchy's special integral formula}

Assume that $f(z)$ does not have poles in the contour and surface enclosed,
where it is assumed to be continuous. At any point $z_{0}$ inside, we have,
by the theorem of residues,%
\begin{equation}
0=\oint \frac{f(z)-f(z_{0})}{z-z_{0}}dx,
\end{equation}%
since $f(z)-f(z_{0})$ is the coefficient of the term of order $-1.$ Assume
further that $f(z_{0})$ is a $2-$form, which is the reason why we put
``special'' in the title. Although $dx$ equals $dz$, we prefer to use $%
f(z)dx $ over $f(z)dz$ to help readers remember this special case. The $%
f(z)dz$ in this paper is not equivalent to the $f(z)dz$ of standard complex
variable calculus. In the next paper, we shall have something equivalent,
but not under the notation $f(z)dz$.

From (30) follows that%
\begin{equation}
\oint \frac{f(z)}{z-z_{0}}dx=\oint \frac{f(z_{0})}{z-z_{0}}dx=\oint \frac{%
f(z_{0})dxdy}{z-z_{0}}dy.
\end{equation}%
Since $f(z_{0})dxdy$ now is a 0-form, we can pull it out of the integral to
further obtain%
\begin{equation}
\oint \frac{f(z)}{z-z_{0}}dx=f(z_{0})dxdy\oint \frac{1}{z-z_{0}}dy^{\prime
}=f(z_{0})dxdy\oint d\phi ^{\prime }=2\pi f(z_{0})dxdy,
\end{equation}%
where we have used that $dy=dy^{\prime }$ (If $f(z_{0})$ were a $0-$form, we
would be dealing with an integral over $\rho $, which is zero). Cauchy's
integral formula then follows:%
\begin{equation}
f(z_{0})=\frac{1}{2\pi dxdy}\oint \frac{f(z)}{z-z_{0}}dx,
\end{equation}

Consider for example the integration%
\begin{equation}
\int_{-\infty }^{+\infty }\text{ }\frac{1}{\left( x^{2}+1\right) }dx=\oint 
\frac{1}{z^{2}+1}=\oint \frac{\frac{1}{z+dxdy}}{z-dxdy}dx
\end{equation}%
on the upper half plane, where $z_{0}=dxdy.$ The evaluation of $1/(z+dxdy)$
at $z=dxdy$ is a $2-$form. The application of (32) is justified. We obtain%
\begin{equation}
\int_{-\infty }^{+\infty }\text{ }\frac{1}{\left( x^{2}+1\right) }dx=2\pi
dxdy\left| \frac{1}{z+dxdy}\right| _{_{z=dxdy}}=\frac{2\pi dxdy}{2dxdy}=\pi .
\end{equation}

In a circle of radius unity centered at $z=0$, the integral 
\begin{equation}
\oint \frac{1}{z(z-\pi )}dz
\end{equation}%
is imaginary. The Cauchy integral formula of this formalism does not apply
here because $1/(z-\pi )$ is real at $z=0$.

Finally consider the integral%
\begin{equation}
\int_{-\infty }^{+\infty }\frac{e^{itx}}{x^{2}+1}dx.
\end{equation}%
We make a closed curve with the upper semicircle at infinity. It contains
the single pole ($0,1$), which implies $z=dxdy$. At $z=dxdy$, the numerator
of%
\begin{equation}
...=\int_{-\infty }^{+\infty }\frac{\frac{e^{tzdxdy}}{z+dxdy}}{z-dxdy}dx.
\end{equation}%
is a 2-form. We apply the theorem and get%
\begin{equation}
2\pi dxdye^{-t}\frac{1}{2dxdy}=\pi e^{-t}.
\end{equation}%
The integrand in (37) is not real, but the integral is.

\subsection{Cauchy's special integral formula for derivatives}

We rewrite Eq. (33), which is valid if $f(z_{0})$ is a $2-$form, as%
\begin{equation}
f(z)=\frac{1}{2\pi dxdy}\oint \frac{f(\zeta )}{\zeta -z}d\chi .
\end{equation}%
Clearly%
\begin{equation}
\frac{\partial f(z)}{\partial x}=\frac{1}{2\pi dxdy}\oint \frac{\partial }{%
\partial x}\frac{f(\zeta )}{\zeta -z}d\chi =\frac{1}{2\pi dxdy}\oint \frac{%
f(\zeta )}{(\zeta -z)^{2}}d\chi
\end{equation}%
and, by successive application,%
\begin{equation}
\frac{\partial ^{n}f(z)}{\partial x^{n}}=\frac{n!}{2\pi dxdy}\oint \frac{%
f(\zeta )}{(\zeta -z)^{n+1}}d\chi ,
\end{equation}%
which we rewrite as%
\begin{equation}
\oint \frac{f(z)}{(z-z_{0})^{n+1}}dx=\frac{2\pi dxdy}{n!}\left( \frac{%
\partial ^{n}f(z)}{\partial x^{n}}\right) _{z=z_{0}}.
\end{equation}

Consider again the integration of (24) around the pole $(0,1)$. In our real
formalism, $z_{0}=dxdy$. 
\begin{equation}
\oint \frac{dx}{\left( z^{2}+1\right) ^{2}}=\oint \frac{\frac{1}{\left(
z+dxdy\right) ^{2}}}{\left( z-dxdy\right) ^{2}}dx.
\end{equation}%
We identify $n=1$ and $f(z)=\left( z+dxdy\right) ^{-2}$ so that%
\begin{equation}
\oint \frac{1}{\left( z^{2}+1\right) ^{2}}dx=2\pi dxdy\left[ \frac{\partial 
}{\partial x}\frac{1}{\left( z+dxdy\right) ^{2}}\right] _{z=dxdy}=\frac{\pi 
}{2}.
\end{equation}

\section{Other theorems}

In the calculus of complex variable, Cauchy's integral formulas are the gate
to the Laurent series, the theorem of residues and the computation of the
residue for poles of order higher than one. We have already covered the
theorem of residues. We assumed at an early point that we would consider
only functions defined by power series (including negative exponents), but
those function will be given to us in abbreviated form, like saying $\sin x$
without stating the series. Occasionally we may want or wish to have the
series explicitly. In the Weierstrass approach that we advocate here, the
perspective is somewhat different from the standard one as we do not need to
use Cauchy's integral formulas.

Given a function $f(x+ydxdy)$ that is defined by a power series (like, say,
a sine function divided by a polynomial which has zero of order $m$ at $%
z_{0} $), it will obviously be given by a power series of the form%
\begin{equation}
f(z)=\frac{a_{-m}}{z^{\prime m}}+\frac{a_{-m+1}}{z^{\prime m-1}}+...+\frac{%
a_{-1}}{z^{\prime }}+a_{0}+a_{1}z^{\prime }+a_{0}z^{\prime 2}+...
\end{equation}%
where $z^{\prime }=z-z_{0}=(x+ydxdy)-(x_{0}+ydxdy).$ The function $z^{\prime
m}f(z)$ is analytic,%
\begin{equation}
z^{\prime m}f(z)=a_{-m}+a_{-m+1}z^{\prime }+...+a_{-1}z^{\prime
m-1}+a_{0}z^{\prime m}+...=\sum_{n=0}^{\infty }b_{n}z^{\prime n}.
\end{equation}%
We now consider $z^{\prime m}f(z)$ as a function in ($x,y$). It is clear that%
\begin{equation}
b_{n}=\frac{1}{n!}\lim_{z\rightarrow z_{0}}\frac{\partial ^{n}[z^{\prime
m}f(z)]}{\partial x^{n}}=a_{n-m}.
\end{equation}%
That is much clearer than the expression for the coefficients of the Laurent
expansion in terms of the integrals in the Cauchy integral formulas. The
Laurent series then takes the form%
\begin{equation}
f(z)=\sum_{n=0}^{\infty }z^{\prime (n-m)}\frac{1}{n!}\lim_{z\rightarrow
z_{0}}\frac{\partial ^{n}[z^{\prime m}f(z)]}{\partial x^{n}}.
\end{equation}

Finally, the formula for the residual when a pole is of order $m$ is a
particular case of (48), namely%
\begin{equation}
a_{-1}=\frac{1}{(m-1)!}\lim_{z\rightarrow z_{0}}\frac{\partial
^{m-1}[z^{\prime m}f(z)]}{\partial x^{m-1}}.
\end{equation}

\section{Concluding remarks}

We have touched only the elementary part of the calculus of complex
variable, which meets the practical interests of most physicists. But it
still should have foundational interest for mathematicians. Those interests
would not have been enough motivation to write this paper. My most important
motivation is what I see as missed opportunity by physicists in overlooking
the worth of the K\"{a}hler calculus for physics.

When endowed with the rich structure conferred by Clifford algebra and K\"{a}%
hler differentiation, differential forms are the words in a unified language
to deal with the requirements of mathematical physics. Their use for the
foundations of quantum mechanics is evident in K\"{a}hler's cited paper (no
negative energy antiparticles, spin without internal space, etc.), which
unfortunately is in German. In order to minimize the need for sufficient
knowledge of that language to follow the arguments through the formulas, it
is helpful to be familiar with the ways of \'{E}. Cartan in dealing with
differential forms, since K\"{a}hler's style is essentially the same. This
style morphed in mid 20th century into the modern one. Though the
recommendation is self-serving, some readers may still wish to consider the
recent book on differential geometry by this author \cite{Book1} in order to
learn to compute with differential forms in Cartan and K\"{a}hler's style.

\end{document}